\documentclass[oneside,final]{csri24}
\usepackage[dvipsnames]{xcolor}

\newcommand{\irm}[1]{{\color{OliveGreen}{IRM:}} {\color{OliveGreen}{#1}}}

\usepackage{tikz}
\usepackage{adjustbox}

\usepackage{mathtools}
\usepackage[top=1.355in,
			bottom=1.355in,
			left=1.5in,
			right=1.5in]{geometry}
\usepackage{amsfonts,
			amsmath,
            amssymb,
			graphicx,
			url}
\usepackage{caption}
\usepackage{subcaption}

\title{Domain decomposition-based coupling of Operator Inference Reduced Order Models via the Schwarz alternating method}

\author{Ian Moore\thanks{Virginia Tech, ianm9123@vt.edu},\and Christopher R. Wentland\thanks{Sandia National Laboratories, crwentl@sandia.gov}\and Anthony Gruber\thanks{Sandia National Laboratories, adgrube@sandia.gov}\and Irina Tezaur\thanks{Sandia National Laboratories, ikalash@sandia.gov}}

\begin{document}

\maketitle


\begin{abstract}
This paper presents and evaluates an approach for coupling subdomain-local reduced order models (ROMs) constructed via non-intrusive operator inference (OpInf) with each other and with subdomain-local full order models (FOMs), following a domain decomposition of the spatial geometry on which a given partial differential equation (PDE) is posed. Joining subdomain-local models is accomplished using the overlapping Schwarz alternating method, a minimally-intrusive multiscale coupling technique that transforms a monolithic problem into a sequence of subdomain-local problems, which communicate through transmission boundary conditions imposed on the subdomain interfaces.  After formulating the overlapping Schwarz alternating method for OpInf ROMs, we evaluate the method's accuracy and efficiency on several test cases involving the heat equation in two spatial dimensions.  We demonstrate that the method is capable of coupling arbitrary combinations of OpInf ROMs and FOMs, and that moderate speed-ups over a monolithic FOM are possible when performing OpInf ROM coupling.
\end{abstract}

\section{Introduction}

Despite advancements in both computer architectures and algorithms, the modeling and simulation of complex physical systems often requires tremendous computational resources. These requirements may preclude many-query analyses such as engineering design or uncertainty quantification.  While projection-based reduced order models (ROMs) have shown promise to mitigate this difficulty, traditional intrusive ROMs, e.g., Galerkin~\cite{Holmes:1996, Sirovich:1987} and least squares Petrov-Galerkin (LSPG) projection ROMs~\cite{Carlberg:2011},  have their own shortcomings, including a lack of systematic refinement mechanisms, a lack of robustness, stability, and accuracy in the predictive regime, and lengthy implementation time requirements.

This paper presents an approach for mitigating the aforementioned difficulties by enabling domain decomposition- (DD-)based coupling of subdomain-local ROMs with each other and/or with subdomain-local full order models (FOMs).  Our approach is based on the following ingredients:  (i) a decomposition of the physical domain of interest into two or more overlapping subdomains, (ii) the construction of subdomain-local ROMs and/or FOMs in each of the subdomains, and (iii) the rigorous coupling of the subdomain-local models via the Schwarz alternating method \cite{Schwarz:1870}.  The Schwarz alternating method is based on the simple idea that, if the solution to a partial differential equation (PDE) is known in two or more regularly shaped domains, these local solutions can be used to iteratively build a solution for the union of the subdomains, with information propagating between the subdomains through carefully constructed transmission boundary conditions (BCs).  We choose the Schwarz alternating method since it has a number of advantages over competing multiscale coupling methods.  These advantages include its concurrent nature, its ability to couple non-conformal meshes with different 
element topologies and different time integrators with different time steps for dynamic problems without introducing non-physical artifacts into the solution, and its non-intrusive implementation into existing codes~\cite{Mota:2017, Mota:2022}. 

Building on our past work in developing the Schwarz alternating method as a means to couple FOMs~\cite{Mota:2017, Mota:2022}, intrusive projection-based ROMs~\cite{Barnett:2022Schwarz} and physics-informed neural networks (PINNs)~\cite{Snyder:2023}, we focus our attention herein on advancing the method to work with a non-intrusive model order reduction (MOR) technique known as operator inference (OpInf)~\cite{ Ghattas_Willcox_2021, kramerOpinf2024,willcox2016opinf}.  Unlike traditional intrusive MOR, which requires access to the underlying FOM code in order to project the governing PDE(s) onto a reduced subspace, OpInf works by assuming a functional form (usually linear or quadratic~\cite{Geelen:2023}) for the ROM in terms of to-be-learned reduced operators, and solving an optimization problem offline for these operators.  This procedure significantly reduces both the development time and the time-to-impact.   

 DD-based FOM-ROM and ROM-ROM couplings such as those proposed herein have
the potential of improving the predictive viability of projection-based ROMs, by enabling the spatial localization of ROMs (via domain decomposition) and the online integration of high-fidelity information into these models (via FOM coupling).
While DD-based couplings between ROMs and FOMs are not new, the majority of the literature on this topic has focused on developing intrusive coupling methods (e.g., Lagrange multipliers, optimization-based coupling) used to couple intrusive ROMs; the interested reader is referred to~\cite{Barnett:2022Schwarz,  Cinquegrana:2011,  deCastro:2023, deCastro:2022,Iollo:2022,Maier:2014, Prusak:2023} and the references therein for more details.  Related past work on data-driven couplings using Schwarz-like methods has focused on intrusive ROMs~\cite{Corigliano:2015, Kerfriden:2013,Kerfriden:2012,Radermacher:2014}, or on utilizing the coupling to accelerate NN training~\cite{LiD3M, LiDeepDDM,Snyder:2023}. The proposed approach is most similar to the recent work by Farcas \textit{et al.}~\cite{Farcas:2023}, which develops a DD-based coupling of subdomain-local OpInf ROMs by learning appropriate reduced operators responsible for the coupling, and demonstrates the method on a formidable three-dimensional (3D) combustion problem.  
In this approach, each subdomain problem is solved once rather than by performing an iteration to convergence as done within our Schwarz framework.  As a result, the subdomain-local solutions must be extended to the full domain and smoothly combined to achieve a continuous solution.  

The remainder of this paper is organized as follows.  In Section~\ref{sec:schwarz}, we describe the overlapping version of the Schwarz alternating method applied to our model problem, the two-dimensional (2D) unsteady heat equation.  In Section~\ref{sec:opinf}, we present some OpInf preliminaries.  
In Section~\ref{sec:software}, we describe our software implementations of the proposed Schwarz-based coupling approach applied to OpInf ROMs, hereafter called the OpInf-Schwarz method, which makes use of the open source {\tt FEniCSx}~\cite{BarattaEtal2023dolfinx} and {\tt OpInf} {\tt Python} libraries.  Numerical results are presented in Section~\ref{sec:results}.  We conclude with a summary and a discussion of future work in Section~\ref{sec:conclusions}.

\section{The Schwarz Alternating Method Applied to Operator Inference ROMs}  \label{sec:schwarz}

In the following sections and subsections, we consider the specific model problem of the 2D heat equation, towards addressing the challenges that are encountered in this novel combination of the Schwarz method and operator inference. We stress that neither of these techniques inherently require an assumption of linearity, noting that the authors of the original operator inference paper~\cite{willcox2016opinf} directly contrasted their method with the linearity assumption of Dynamic Mode Decomposition (DMD) \cite{DMD}. While non-linearity will likely present further challenges to overcome for the OpInf-Schwarz coupling method, the heat equation provides an initial unsteady test problem that can suggest the feasibility of the method for more complicated problems. 

\subsection{Schwarz Alternating Method Preliminaries}
 \label{sec:schwarz_details}

\begin{figure}[htbp!]
        \begin{center}{\includegraphics[width=0.4\textwidth]
    {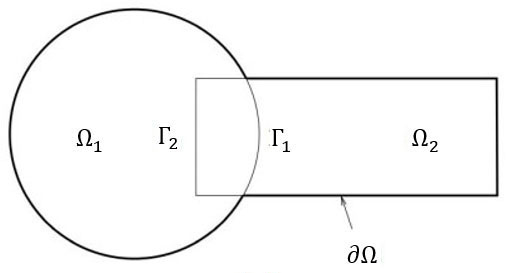}}
        \end{center}
        \caption{Illustration showing an overlapping domain
    decomposition of a 2D domain $\Omega$ for the application of the Schwarz alternating method. Note $\Gamma_1, \Gamma_2 \not\subset\partial\Omega$.}
        \label{fig:dd}
\end{figure}

Consider the heat equation specified as, 
\begin{equation} \label{eq:heat_pde}
\begin{array}{rcll}
 \dot{u}(x,t) - \Delta u(x,t) &=& 0, & \text{in } \Omega \times [0, T], \\
    u(x,t) &=& g(x), & \text{on } \partial \Omega \times [0,T],\\
    u(x,0) & = & v(x), & \text{in } \Omega,
    \end{array}
\end{equation}
where $\Omega \in \mathbb{R}^p$ is an open bounded domain for $p=1,2,3$ with boundary $\partial \Omega$,  $g(x)$ are is a given boundary condition function, $v(x)$ defines the initial condition, and $T>0$. Suppose we decompose the domain $\Omega$ into two overlapping subdomains $\Omega_1$ and $\Omega_2$, such that $\Omega_1 \cup \Omega_2 = \Omega$ and $\Omega_1 \cap \Omega_2 \neq \emptyset$ and, for later use, define $\overline{\Omega_1}$ as the closure of $\Omega_1$ and similarly for $\Omega_2$.  Suppose also that we decompose the time interval $[0,T]$ into a set of non-overlapping time intervals $I_n = [t_n, t_{n+1}]$, where $T \ge t_{n+1} > t_n \ge 0$,
so that $\cup_{n}I_n = [0,T]$.  For a given time interval $I_n$, the overlapping Schwarz algorithm solves the following sequence of subdomain-local problems: 
\begin{equation} \label{eq:generic_schwarz_iter_Omega1}
    \left \{
    \begin{array}{rcll}
    \dot{u}_1^{(k+1)} - \Delta u_1^{(k+1)} &=& 0, & \text{in } \Omega_1 \times [t_n, t_{n+1}] \\ 
     u_1^{(k+1)} &=& g, & \text{on } (\partial \Omega \cap \overline{\Omega_1}) \times [t_n,t_{n+1}]\\
     u_1^{(k+1)} &=&  u_2^{(k)},  &\text{on } \Gamma_1 \times [t_n,t_{n+1}], 
    \end{array}
    \right.
\end{equation}
and 
\begin{equation} \label{eq:generic_schwarz_iter_Omega2}
    \left \{
    \begin{array}{rcll}
    \dot{u}_2^{(k+1)} - \Delta u_2^{(k+1)} &=& 0, & \text{in } \Omega_2 \times [t_n, t_{n+1}] \\ 
     u_2^{(k+1)} &=& g, & \text{on } (\partial \Omega \cap \overline{\Omega_2}) \times [t_n,t_{n+1}]\\
     u_2^{(k+1)} &=& u_1^{(k+1)},  &\text{on } \Gamma_2 \times [t_n,t_{n+1}], 
    \end{array}
    \right.
\end{equation}
for Schwarz iteration $k = 0, 1, 2, ...,$ subject to  initial conditions $u_1(x,0) = v|_{\Omega_1}$ and $u_2(x,0) = v|_{\Omega_2}$. In~\eqref{eq:generic_schwarz_iter_Omega1}, $u_i$ for $i=1,2$ denotes the solution in subdomain $\Omega_i$, and $\Gamma_i$ is the so-called Schwarz boundary (see Figure~\ref{fig:dd}).  It is common to set $u_1^{(0)} = v|_{\partial \Omega_1}$ on $\Gamma_1$ and $u_2^{(0)} = v|_{\partial \Omega_2}$ on $\Gamma_2$ when initializing the Schwarz iteration process, to ensure solution compatibility with the initial condition. The iterative process in~\eqref{eq:generic_schwarz_iter_Omega1} and~\eqref{eq:generic_schwarz_iter_Omega2} continues until a set of pre-determined criteria are met.  In the present work, convergence criteria are based on the Euclidean norm of the solution differences between consecutive Schwarz iterations; that is, Schwarz is deemed converged when $\epsilon_{\text{abs}}^{(k)} < \delta_{\text{abs}}$ and $\epsilon_{\text{rel}}^{(k)} < \delta_{\text{rel}}$ for some pre-specified Schwarz tolerances $\delta_{\text{abs}}, \delta_{\text{rel}} > 0$, where 
\begin{equation} \label{eq:conv_criterion_abs}
    \epsilon_{\text{abs}}^{(k)} := \sqrt{|| u_1^{(k)} - u_1^{(k-1)}||^2 + || u_2^{(k)} - u_2^{(k-1)}||^2},
\end{equation}
and 
\begin{equation} \label{eq:conv_criterion_rel}
    \epsilon_{\text{rel}}^{(k)} := \sqrt{\frac{|| u_1^{(k)} - u_1^{(k-1)}||^2}{||u_1^{(k)}||^2} + \frac{|| u_2^{(k)} - u_2^{(k-1)}||^2}{||u_2^{(k)}||^2}},
\end{equation}
for Schwarz iteration $k=1,2,...$, until convergence.

A key advantage of the Schwarz alternating method is that it allows for the subdomains $\Omega_1$ and $\Omega_2$ to be discretized using different meshes and/or element types~\cite{Mota:2017, Mota:2022}.  In the case where $\Omega_1$ and $\Omega_2$ are discretized using different meshes and do not have a coincident interface, applying the Schwarz boundary condition on the Schwarz boundaries $\Gamma_i$ requires the construction of a projection operator; this can be done via a simple application of finite element interpolation functions readily available in most codes~\cite{ Mota:2023,Mota:2022}. 

The Schwarz iteration process~\eqref{eq:generic_schwarz_iter_Omega1}--\eqref{eq:generic_schwarz_iter_Omega2} is converged within each time interval $[t_n, t_{n+1}]$ before moving on to the next time interval.  A key advantage of this time-stepping approach is that different time integrators and time steps can be used in different subdomains; for a detailed discussion of Schwarz time-stepping and related machinery, the reader is referred to~\cite{Mota:2022, Mota:2023}. In the numerical experiments presented herein, we restrict attention to the case where all subdomains have the same time integrator and time step, as our goal is to assess the method's viability when coupling subdomain-local operator inference-based ROMs with each other and with subdomain-local FOMs while eliminating other confounding factors.

\subsubsection{FOM-FOM Schwarz Coupling}
\label{sec:FOM-FOM_schwarz}
Continuing the example of the heat equation from~\eqref{eq:heat_pde}, a spatially discretized monolithic FOM for the heat equation typically appears in the following form after a boundary lift:  
\begin{equation}
\label{eq:heat_eq_discretized}
    \dot{\mathbf{x}} = \mathbf{K} \mathbf{x} + \mathbf{B} \mathbf{g},
\end{equation}
where $\mathbf{x} \in \mathbb{R}^{N}$ is a discretized vector corresponding to the unconstrained state degrees of freedom (DoFs), and $\mathbf{g} \in \mathbb{R}^{m}$ discretizes the Dirichlet boundary condition. The matrix $\mathbf{K} \in \mathbb{R}^{N\times N}$ comes from discretizing the continuous Laplace operator $\Delta$, and $\mathbf{B} \in \mathbb{R}^{N \times m}$ deals with effects of the boundary condition. The full state representation $\mathbf{u} \in \mathbb{R}^{N + m}$ for all DoFs is obtained by augmenting the unconstrained solution $\mathbf{x}$ with the known boundary condition $\mathbf{g}$.

If a domain decomposition such as in Figure \ref{fig:dd} is imposed, we may consider the subdomain-local discretized problems: 
\begin{equation}
\label{eq:heat_eq_matrices_DD}
    \dot{\mathbf{x}}_i = \mathbf{K}_i \mathbf{x}_i + \mathbf{B}_i \mathbf{y}_i,   \;\;\; i = 1,2,
\end{equation}
on $\Omega_1$ and $\Omega_2$. Let $N_i$ and $m_i$ denote the number of finite element (FE) nodes on the interior and boundary of $\Omega_i$, respectively. Note that the boundary of $\Omega_i$ includes $\Gamma_i$, as shown in Figure \ref{fig:dd}. We introduce the vector $\mathbf{y}_i \in \mathbb{R}^{m_i}$ to handle the subdomain-local boundary by defining $\mathbf{y}_i = [\mathbf{g}_i \;\mathbf{\gamma}_i]^T$, where $\mathbf{g}_i$ discretizes $g(x)$ on $(\partial \Omega \cap \overline{\Omega_i})$ and $\gamma_i$ discretizes the subdomain-local boundary condition on $\Gamma_i$.

The remaining quantities in~\eqref{eq:heat_eq_matrices_DD} follow naturally from the monolithic terms in \eqref{eq:heat_eq_discretized}, with subdomain-local state representation $\mathbf{x}_i \in \mathbb{R}^{N_i}$ as well as stiffness and boundary matrices $\mathbf{K}_i \in \mathbb{R}^{N_i\times N_i}$  and $\mathbf{B}_i \in \mathbb{R}^{N_i \times m_i}$. We now use the Schwarz method to set up the subdomain-local problems as follows: 
\begin{equation} \label{eq:FOM_schwarz_iter_Omega1}
    \left \{
    \begin{array}{rcll}
    \dot{\mathbf{x}}_1^{(k+1)} &=& \mathbf{K}_1 \mathbf{x}_1^{(k+1)} + \mathbf{B}_1 \mathbf{y}_1^{(k+1)} \\
     \gamma_1^{(k+1)} &=&  \mathbf{x}_2^{(k)} \big|_{\Gamma_1},   
    \end{array}
    \right.
\end{equation}
and 
\begin{equation} \label{eq:FOM_schwarz_iter_Omega2}
    \left \{
    \begin{array}{rcll}
    \dot{\mathbf{x}}_2^{(k+1)} &=& \mathbf{K}_2 \mathbf{x}_2^{(k+1)} + \mathbf{B}_2 \mathbf{y}_2^{(k+1)} \\
     \gamma_2^{(k+1)} &=&  \mathbf{x}_1^{(k+1)} \big|_{\Gamma_2},   
    \end{array}
    \right.
\end{equation}
for Schwarz iteration $k = 0, 1, 2, ...,$ until convergence following \eqref{eq:conv_criterion_abs}-\eqref{eq:conv_criterion_rel}. For work presented in this paper, we solve equations \eqref{eq:FOM_schwarz_iter_Omega1}-\eqref{eq:FOM_schwarz_iter_Omega2} using the finite element method in space and the finite difference method in time.  We take the initial condition to be the interpolation of $u(x,0) = v(x)$ into the respective finite element spaces. This specific formulation is valid for conformal spatial meshes only.

\subsection{Operator Inference Preliminaries} \label{sec:opinf}

Operator Inference is a data-driven, non-intrusive, projection-based method for model order reduction, developed by Peherstorfer and Willcox~\cite{willcox2016opinf} as an alternative to standard Galerkin projection ROMs. Similarly to other Galerkin projection ROMs, the method begins by constructing a reduced basis from data. OpInf then diverges from the standard method in the construction of the reduced operators, which are estimated using data through a regression problem. 

\subsubsection{Proper Orthogonal Decomposition}

To construct a reduced basis, we use proper orthogonal decomposition (POD)~\cite{Holmes:1996, Kunisch2001POD,Sirovich:1987}, though one might also use other methods, e.g., the reduced basis method with a greedy algorithm~\cite{RozzaGreedy2014}. To perform POD, we require snapshots from some data source or FOM.   

Continuing with the heat equation example, after the monolithic problem~\eqref{eq:heat_eq_discretized} has been solved in time for $\tau$ separate states inclusively between time $t = 0$ and $t = T$ (i.e., $0 = t_1 < t_2 < \dots < t_\tau = T$), we obtain a collection of unconstrained state snapshots, 
\begin{equation}
\mathbf{X} = [\mathbf{x}(t_1), \mathbf{x}(t_2), \dots, \mathbf{x}(t_\tau)] \in \mathbb{R}^{N 
\times \tau},
 \label{eq:snapshots}
\end{equation}
with rank$(\mathbf{X}) = d > 0$. $\mathbf{X}$ yields a singular value decomposition (SVD) $\mathbf{X} = \mathbf{\Psi} \mathbf{\Sigma} \mathbf{\Phi}^*$,  where $\mathbf{\Psi} \in \mathbb{R}^{N \times N}$ is an orthogonal matrix whose first $d$ columns form a basis for the column space of $\mathbf{X}$. The columns of $\mathbf{\Psi}$, $\mathbf{\psi}_i$ for $i = 1, \dots r$ with $r \leq d$ form an optimal $r$-dimensional basis for the columns of $\mathbf{X}$ in an $\ell^2$ sense. The restriction of $\mathbf{\Psi}$ to its first $r$ columns is referred to as $\mathbf{\Psi}_r \in \mathbb{R}^{N\times r}$ going forwards.

\subsubsection{Operator Inference}
\label{sec:Opinf_Details}
Operator inference is a non-intrusive, data-driven, project- ion-based model order reduction technique which learns low-dimensional operators that can be used to approximate the output of a given FOM, which, in this paper, is the heat equation described in~\eqref{eq:heat_eq_discretized}. In a classical intrusive Galerkin ROM, the operator $\mathbf{K}_r \in \mathbb{R}^{r\times r}$ is the operator which best represents the action of $\mathbf{K}$ in the reduced space of dimension $r$ defined by the basis $\mathbf{\Psi}_r$. That is,
\begin{equation}
 \mathbf{K}_r = \mathbf{\Psi}_r^T \mathbf{K}  \mathbf{\Psi}_r.
 \label{eq:Galerkin_ROM_Operator}
\end{equation}
This is intrusive because it requires access to the FOM matrix $\mathbf{K}$, and, practically speaking, the code that produced $\mathbf{K}$.  

The key difference between operator inference
and traditional intrusive projection-based MOR is that, in the former approach,
the reduced operators are not created via intrusive projection, but inferred directly using available snapshot data.  OpInf is based on the observation that a projection-based ROM derived from a FOM with polynomial nonlinearities will possess the same algebraic structure as the FOM.  

For the problem considered herein, the 2D heat equation, it is straightforward to see that a projection-based ROM preserves the linear algebraic structure of the FOM,~\eqref{eq:heat_eq_discretized}. 
With this observation, we now set up the semi-discretized monolithic OpInf ROM problem:
\begin{equation} \label{eq:opinf_rom}
    \dot{\hat{\mathbf{x}}} = \hat{\mathbf{K}} \hat{\mathbf{x}} + \hat{\mathbf{B}} \mathbf{g}.
\end{equation}
In order to obtain the reduced operators without access to FOM matrices, one solves a regression problem for the 
reduced matrices
$\hat{\mathbf{K}} \in \mathbb{R}^{r \times r}$ and $\hat{\mathbf B} \in \mathbb{R}^{r \times m}$. In particular, given $j \leq \tau$ steps of FOM training data $\mathbf{x}_p = \mathbf{x}(t_p)$, $p = 1, \dots, j$, minimize,
\begin{equation}
\label{eq:OpInf_2_norm}
    \min_{\hat{\mathbf{K}}, \hat{\mathbf{B}}} \sum_{p=1}^j \| \dot{\hat{\mathbf{x}}}_p - \hat{\mathbf{K}} \hat{\mathbf{x}}_p - \hat{\mathbf{B}} \mathbf{g}\|_2^2.
\end{equation}
Within the training data, we define the ROM representation of the state at time $t_p$ to be $\hat{\mathbf{x}}_p = \mathbf{\Psi_r^T} \mathbf{x}_p \in \mathbb{R}^r$. The time derivative of the ROM state, $\dot{\hat{\mathbf{x}}}_{p} \in \mathbb{R}^r$, may not be explicitly available and can instead be estimated from available data using a difference method. The boundary information, $\mathbf{g} \in \mathbb{R}^m$, is the same as in the monolithic FOM. 

We now briefly contextualize this in a classical Galerkin ROM context. A typical projection-based ROM creates matrices which are in some way representative of those constructed for the FOM, but OpInf does not start with the assumption that these FOM matrices are available. A more natural comparison is between OpInf and the standard Galerkin ROM for the choice of a particular shared basis $\mathbf{\Psi}_r$. 

It is guaranteed that, under the assumptions that the time-stepping scheme is convergent as $dt \to 0$, the approximation $\dot{\hat{\mathbf{x}}}_j$ converges to the true time derivative at time $t_j$ as $dt \to 0$, and linear independence of supplied data ($d \geq \tau$, the total number of snapshots), then, for all $0 < \epsilon \in \mathbb{R}$ there exists a timestep $dt$ for the FOM and a $r \leq d$ such that,
\begin{equation*}
    \| \mathbf{K}_r - \hat{\mathbf{K}} \|_F < \epsilon.
\end{equation*}
The same result holds for $\hat{\mathbf{B}}$~\cite{willcox2016opinf}. In the case that a technique known as re-projection is used for the fully discretized problem, it can be shown that the learned OpInf ROM system recovers
the standard Galerkin projected intrusive ROM~\cite{Peherstorfer:2020}.


While this overview has focused on the heat equation, we again stress that operator inference is general to low-order polynomial non-linearity, with results previously shown for a polynomial non-linearity of degree three in a one-dimensional (1D) nuclear reactor model in~\cite{willcox2016opinf}, for a 2D single injector combustion model in~\cite{McQuarrie2021combustion}, and for a 3D rotating detonation rocket engine in~\cite{Farcas:2023}, among others.

\subsection{Opinf-Schwarz Method: FOM-ROM Coupling}
\label{sec:FOM-ROM_Schwarz}
We are now in a position to state the Opinf-Schwarz method for the heat equation specified in \eqref{eq:heat_pde}. We state the problem formulation for the case of two overlapping subdomains $\Omega_1$ and $\Omega_2$ with $\Omega_1 \cup \Omega_2 = \Omega$ as in Figure \ref{fig:dd} where we place a FOM on $\Omega_1$, an OpInf ROM on $\Omega_2$, and couple the problems via the Schwarz method. The Opinf-Schwarz formulation can be extended easily to arbitrary numbers of subdomains and combinations of model couplings. For simplicity, we assume that we have access to monolithic FOM training data across the entire domain $\Omega$. The dimensional and subdomain notation presented here follows from that defined for the FOM-FOM coupling in Section \ref{sec:FOM-FOM_schwarz}.

The OpInf problem on $\Omega_2$ requires some details. Given $j \leq \tau$ steps of monolithic data $\mathbf{u}(t_p) \in \mathbb{R}^{N + m}, p = 1,\dots,j$ extract the unconstrained state representation on $\Omega_2$, $\mathbf{x}_2(t_p) \in \mathbb{R}^{N_2}$ as well as the the vector $\gamma_2(t_p)$, which is the state information along the Schwarz boundary $\Gamma_2$ at time $t_p$. The subdomain local boundary on $\Omega_2$ is completed by defining the vector $\mathbf{y}_2(t_p) = [\mathbf{g}_2 \; \gamma_2(t_p)]^T \in \mathbb{R}^{m_2}$. We carry out POD on the collection of unconstrained snapshots on $\Omega_2$, $\mathbf{x}_2(t_p), p = 1, \dots, j$ to obtain the $\Omega_2$-local basis $\mathbf{\Psi}_{r,2} \in \mathbb{R}^{N_2 \times r}$ for a choice of basis dimension $r$. Once the basis is obtained, we solve the regression problem, 
\begin{equation}
\label{eq:opinf-schwarz_regression_problem}
    \min_{\hat{\mathbf{K}}_2, \hat{\mathbf{B}}_2} \sum_{p=1}^j \| \dot{\hat{\mathbf{x}}}_2(t_p) - \hat{\mathbf{K}}_2 \hat{\mathbf{x}}_2(t_p) - \hat{\mathbf{B}}_2 \mathbf{y}_2(t_p)\|_2^2,
\end{equation} 
for the reduced operators $\hat{\mathbf{K}}_2 \in \mathbb{R}^{r \times r}$ and $\hat{\mathbf{B}}_2 \in \mathbb{R}^{r \times m_i}$ as in Section \ref{sec:opinf}.

After setting up the FOM on $\Omega_1$ exactly as in Section \ref{sec:FOM-FOM_schwarz},  the problems may be coupled using the Schwarz method as follows:

\begin{equation} \label{eq:FOM-ROM_schwarz_iter_Omega1}
    \left \{
    \begin{array}{rcll}
    \dot{\mathbf{x}}_1^{(k+1)} &=& \mathbf{K}_1 \mathbf{x}_1^{(k+1)} + \mathbf{B}_1 \mathbf{y}_1^{(k+1)} \\
     \gamma_1^{(k+1)} &=&  \mathbf{\Psi}_{r,2}\hat{\mathbf{x}}_2^{(k)} \big|_{\Gamma_1},   
    \end{array}
    \right.
\end{equation}
and 
\begin{equation} \label{eq:FOM-ROM_schwarz_iter_Omega2}
    \left \{
    \begin{array}{rcll}
    \dot{\hat{\mathbf{x}}}_2^{(k+1)} &=& \hat{\mathbf{K}}_2 \hat{\mathbf{x}}_2^{(k+1)} + \hat{\mathbf{B}}_2 \mathbf{y}_2^{(k+1)} \\
     \gamma_2^{(k+1)} &=&  \mathbf{x}_1^{(k+1)} \big|_{\Gamma_2},   
    \end{array}
    \right.
\end{equation}
for Schwarz iteration $k = 0,1,2, \dots$ until convergence following \eqref{eq:conv_criterion_abs}-\eqref{eq:conv_criterion_rel}. $\mathbf{\Psi}_{r,2}\hat{\mathbf{x}}_2^{(k)}$ is the reconstruction of the ROM state $\hat{\mathbf{x}}_2^{(k)}$ on the FOM mesh originating from the data source. 

The specific formulation described above is valid for conformal meshes only; some interpolation, evaluation or projection scheme would need to be employed for non-conformal meshes depending on the choice of spatial discretization for non-conformal meshes. Creating a ROM-ROM coupling in this setting only requires finding the reduced operators on $\Omega_1$ and reconstructing the values found on the Schwarz boundary onto the FOM mesh. We extend this approach to arbitrary numbers of subdomains by traversing the domains 
in a round-robin fashion and use the most up-to-date Schwarz boundary information available at every Schwarz iteration. 



\section{Software Implementations} \label{sec:software}
We now discuss specific choices made in the implementation of the Opinf-Schwarz method for this paper, in particular those which are not essential to the method itself. Our results include Schwarz coupling implementations for the FOM-FOM (Section \ref{sec:FOM-FOM_schwarz}) FOM-ROM (Section \ref{sec:FOM-ROM_Schwarz}) and ROM-ROM (Adaptation of Section \ref{sec:FOM-ROM_Schwarz}) scenarios. In our implementation, the generic FOM now specifically refers to a finite element simulation, and ROM refers to operator inference. All ordinary differential equations (ODEs), such as those presented for the Opinf-Schwarz method in equations \eqref{eq:FOM-ROM_schwarz_iter_Omega1}-\eqref{eq:FOM-ROM_schwarz_iter_Omega2}, are discretized in time using backward Euler. 

We use the open source FE library {\tt FEniCSx} in {\tt Python} for all FE simulations. {\tt FEniCSx} is a FE problem solving environment with the ability to define variational forms close to actual mathematical notation~\cite{AlnaesEtal2014}, along with various tools~\cite{BarattaEtal2023dolfinx,BasixJoss,ScroggsEtal2022} to enable the creation of a FE space on which computations can be performed. Boundary conditions are enforced strongly through a lifting method leading to FOM systems similar to~\eqref{eq:heat_eq_discretized}. 

Operator inference is performed with the {\tt Python} library {\tt OpInf}, a  {\tt PyPI Python} package available at https://pypi.org/project/opinf/. {\tt OpInf} is a set of tools to facilitate learning polynomial reduced order models, applied in this paper as described in Section~\ref{sec:opinf}. The required data is drawn from a monolithic FE simulation on $\Omega$, also performed in 
{\tt{FEniCSx}}. When solving the OpInf regression problem~\eqref{eq:opinf-schwarz_regression_problem}, we estimate the required state derivative information $\dot{\hat{\mathbf{x}}}_2$ for  through a first order backward difference of available state data.

The numerical stability of OpInf models is a significant concern and area of active research. OpInf models are not guaranteed to be stable, but 
there are a variety of strategies that can help improve stability of reduced models -- independently of OpInf, data pre-processing is commonly used in ROMs to create a model for a transformed set of snapshots. In this paper, we use a centering approach where we build a model for snapshot data $\mathbf{x} - \mathbf{\bar{x}}$, the mean-centered snapshots. 

Specific to OpInf, we also employ a regression regularization strategy for equation~\eqref{eq:OpInf_2_norm}. When solving an OpInf regression problem like~\eqref{eq:OpInf_2_norm}, we add a regularization form which converts the problem into the following form: 
\begin{equation}
\label{eq:OpInf_regularization}
    \min_{\hat{\mathbf{K}}, \hat{\mathbf{B}}} \sum_{p=1}^j \| \dot{\hat{\mathbf{x}}}_p - \hat{\mathbf{K}} \hat{\mathbf{x}}_p - \hat{\mathbf{B}} \mathbf{g}\|_2^2 + \lambda^2 \left(\| \hat{\mathbf{K}}\|^2_F + \|\hat{\mathbf{B}} \|^2_F\right).
\end{equation}
By adding this regularization, we penalize the large entries of our ROM matrices $\hat{\mathbf{K}}, \hat{\mathbf{B}}$. In our testing, this centering and regularization strategy has been sufficient to produce a stable ROM ODE. For our results in Section~\ref{sec:results}, all tables and figures have been produced with regularization parameter $\lambda = 10^{-2}$. Regularization is a common necessity in operator inference, see for example~\cite{McQuarrie2021combustion,Sawant:2023} for related approaches to regularization of OpInf models. 

Finally, the Opinf-Schwarz method that couples these problems has been implemented by the authors of this paper in accordance with Section ~\ref{sec:FOM-FOM_schwarz} and Section~\ref{sec:FOM-ROM_Schwarz}. To simplify the boundary transmission of the Schwarz process, only conformal meshes are used for the differing domains. The output of the Schwarz methods (FOM-FOM, ROM-ROM) as presented in Section \ref{sec:schwarz} are unconstrained state vectors on the interiors of $\Omega_1$ and $\Omega_2$, when what is desired is a single approximation of the monolithic state across $\Omega$. We have made the choice to reconcile the solutions in the following manner. 

Taking the FOM-FOM coupling from Section \ref{sec:FOM-FOM_schwarz} as an example, the Opinf-Schwarz method produces $\{\mathbf{x}_1(t_p)\}_{p=1}^{\tau}$ and $\{\mathbf{x}_2(t_p)\}_{p=1}^{\tau}$. We define the merged solution vector $\mathbf{u}(t_p)_{\text{SCHWARZ}} \in \mathbb{R}^N$ to be equal to the DoFs of $\mathbf{x}_i(t_p)$ on those entries associated with the interior of $\Omega_i \setminus \Omega_1 \cap \Omega_2$ and equal to the known boundary condition $\mathbf{g}$ on entries associated with $\partial \Omega$. For those entries associated with the overlap region $\Omega_1 \cap \Omega_2$, we assign the value of either the entry belonging to $\mathbf{x}_1(t_p)$ or $\mathbf{x}_2(t_p)$ depending on whether the physical node for that entry is closer to the center of $\Omega_1$ or $\Omega_2$.

\section{Results} \label{sec:results}

For the purposes of an initial test of OpInf-Schwarz, we are most interested if the model matches the standard behavior of ROMs and the Schwarz method. We expect that, as the ROM dimension $r$ and the data allowance increases, the error should approach 0, and that as we increase the size of the Schwarz overlap, the number of Schwarz iterations to converge per timestep will decrease. Lastly, as for any ROM, we expect a significant speed-up in computation time for the price we pay in accuracy. 

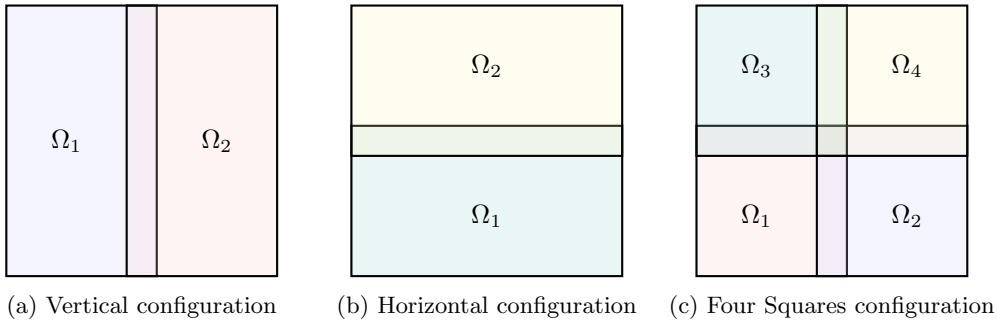
\begin{figure}[h]
         \def\squareSize{2}
         \def\overlap{0.4}

     \begin{subfigure}[c]{0.32\textwidth}
         \centering
         \begin{tikzpicture}


         \pgfmathsetmacro{\adjustedSize}{\squareSize - \overlap}

         \draw[fill=blue!20, fill opacity=0.2, draw=black, thick] (0, 0) rectangle (\squareSize, {2*\squareSize-\overlap});
         \draw[fill=red!20, fill opacity=0.2, draw=black, thick] ({\squareSize-\overlap}, 0) rectangle ({2*\squareSize-\overlap}, {2*\squareSize-\overlap});

         \node at ({\adjustedSize/2}, {\squareSize - \overlap/2}) {$\Omega_1$};
         \node at ({\squareSize+\adjustedSize/2}, {\squareSize - \overlap/2}) {$\Omega_2$};

         \end{tikzpicture}
         \caption{Vertical configuration}
         \label{schematic:vertical}
     \end{subfigure}
     \begin{subfigure}[c]{0.32\textwidth}
         \centering
         \begin{tikzpicture}


         \pgfmathsetmacro{\adjustedSize}{\squareSize - \overlap}

         \draw[fill=Emerald!40, fill opacity=0.2, draw=black, thick] (0, 0) rectangle ({2*\squareSize-\overlap}, \squareSize);
         \draw[fill=Goldenrod!40, fill opacity=0.2, draw=black, thick] (0, {\squareSize-\overlap}) rectangle ({2*\squareSize-\overlap}, {2*\squareSize-\overlap});

         \node at ({\squareSize - \overlap/2}, {\adjustedSize/2}) {$\Omega_1$};
         \node at ({\squareSize - \overlap/2}, {\squareSize + \adjustedSize/2}) {$\Omega_2$};

         \end{tikzpicture}
         \caption{Horizontal configuration}
         \label{schematic:horizontal}
     \end{subfigure}
     \begin{subfigure}[c]{0.32\textwidth}
         \centering
         \begin{tikzpicture}

         \def\squareplus{\squareSize + \overlap}
         \pgfmathsetmacro{\squareminus}{\squareSize - \overlap}

         \draw[fill=red!20, fill opacity=0.2,draw=black, thick] (0, 0) rectangle (\squareSize, \squareSize);
         \draw[fill=blue!20, fill opacity=0.2,draw=black, thick] ({\squareSize-\overlap}, 0) rectangle ({2*\squareSize-\overlap}, \squareSize);
         \draw[fill=Emerald!40, fill opacity=0.2,draw=black, thick] (0, {\squareSize-\overlap}) rectangle (\squareSize, {2*\squareSize-\overlap});
         \draw[fill=Goldenrod!40, fill opacity=0.2,draw=black, thick] ({\squareSize-\overlap}, {\squareSize-\overlap}) rectangle ({2*\squareSize-\overlap}, {2*\squareSize-\overlap});

         \node at ({\squareminus/2}, {\squareminus/2}) {$\Omega_1$};
         \node at ({\squareSize + \squareminus/2}, {\squareminus/2}) {$\Omega_2$};
         \node at ({\squareminus/2}, {\squareSize + \squareminus/2}) {$\Omega_3$};
         \node at ({\squareSize + \squareminus/2}, {\squareSize + \squareminus/2}) {$\Omega_4$};

         \end{tikzpicture}
         \caption{Four Squares configuration}
         \label{schematic:quadrants}
     \end{subfigure}
        \caption{Overlapping domain decomposition configurations }
        \label{fig:schematics}
\end{figure}
We investigate each of these points in several different domain decomposition scenarios as shown to Figure~\ref{fig:schematics}. For a global domain, we choose the square $\Omega = [-1,1] \times [-1,1] \in \mathbb{R}^2$ on which a monolithic FE simulation is carried out. The monolithic boundaries are referred to as $\partial\Omega_L$ for the left, $\partial\Omega_R$ for the right, and  $\partial\Omega_T$ and  $\partial\Omega_B$ for the top and bottom of the domain respectively. We divide the domain into either two rectangles or four squares, each with an equal Schwarz overlap measured in the number of elements overlapping each neighboring subdomain. For these Schwarz overlaps, explicit communication is only between neighbors which share at least one complete edge --hence, in the Four Squares configuration, boundary information can be communicated to a horizontal or vertical neighbor, but not a diagonal neighbor (e.g. $\Omega_1$ communicates directly with $\Omega_2$ and $\Omega_3$ but not $\Omega_4$ in Figure \ref{schematic:quadrants}). The Schwarz subdomains are updated in ascending numerical order.  

The specific problem we test the Opinf-Schwarz method on is the heat equation as specified in~\eqref{eq:heat_pde}, so we fix the following quantities for the sake of comparison. Our time domain is $[0, 1]$, partitioned with a constant time step of $dt = 0.01$ leading to 101 steps, inclusively, from the initial condition at $t_1 = 0$ to the final evaluation at $t_T = 1$. The monolithic domain is be split into a $50\times50$ grid, and then triangular Lagrangian elements of degree 1 is formed by splitting the grid diagonally. These choices lead to 2601 FE nodes on the monolithic mesh, inclusive of the boundary. Within the Schwarz method, our relative~\eqref{eq:conv_criterion_rel} and absolute~\eqref{eq:conv_criterion_abs} tolerances for convergence are both be fixed at $10^{-10}$. 


\paragraph{Errors} Our typical measure of error for OpInf-Schwarz is the average pointwise relative error, i.e.,
\[
\mathbf{E}_{\ell^2}^{\text{avg}} = \frac{1}{\tau-1}\sum_{p = 2}^\tau \frac{ \|\mathbf{u}(t_p)_{ \text{MONO}} - \mathbf{u}(t_p)_{ \text{SCHWARZ}} \|_2}{ \| \mathbf{u}(t_p)_{ \text{MONO}} \|_2},
\]
As in Section \ref{sec:software}, $\mathbf{u}(t_p)_{ \text{SCHWARZ}}$ refers to the solution via the Schwarz method merged from the component subdomains at time $t_p$, and  $\mathbf{u}(t_p)_{ \text{MONO}}$ refers to the monolithic finite element  simulation on $\Omega$ at time $t_p$. 
The error determination excludes the given initial condition, but includes the boundary condition for computational convenience.

The projection error of the snapshots onto the POD basis in subdomain $\Omega_i$ (dependent on basis dimension $r$) is here defined as, 
\[
\mathbf{E}_{\text{proj}}^{\Omega_i} = \frac{\| \mathbf{X}_i - \mathbf{\Psi}_{r,i} \mathbf{\Psi}_{r,i}^T \mathbf{X}_i \|_F}{\|\mathbf{X}_i\|_F},
\]
where $\mathbf{X}_j$ (chosen to recall the unconstrained state variable $\mathbf{x}$, which the POD basis is built for) are the snapshots associated with all solved for DoFs in $\Omega_i$ over the entire time domain, and $\mathbf{\Psi}_{r,i}$ is the $r$-dimensional basis associated with the solved for DoFs of $\Omega_i$. This gives a measure of the approximation quality of the POD basis in each subdomain. When averaged across all subdomains, this is reported using the symbol,  $\mathbf{E}_{\text{proj}}^{\text{avg}}$.

Lastly, we also use the maximum error in Section~\ref{sec:results_time_var}, defined as 
\[
\mathbf{E}_{\ell^2}^{\text{max}} = \max_p \frac{ \|\mathbf{u}(t_p)_{ \text{MONO}} - \mathbf{u}(t_p)_{ \text{SCHWARZ}} \|_2}{ \| \mathbf{u}(t_p)_{ \text{MONO}} \|_2}.
\]
\paragraph{Quantities Reported in Tables}
In addition to errors, we report the following quantities. 
\begin{enumerate}
    \item ``Avg S.I.'' is the mean number of Schwarz iterations required to reach convergence across all time steps.
    \item ``Overlap'' is the amount by which the $50 \times 50$ grid is shared by neighboring subdomains before being broken up into triangular elements. For example, in the Four Squares configuration of Figure \ref{schematic:quadrants}, if Overlap = 5, then the subdomain-specific grid for $\Omega_1$ shares five rows of the original grid with $\Omega_3$ and five columns with $\Omega_2$. The overlap with $\Omega_4$ is incidental. 
    \item ``Time'' is the average wall clock time taken to run the online model on identical hardware, measured in seconds. The number of runs averaged are mentioned in each table. This does not include any time spent setting up matrices, generating data, or similar offline tasks.
    \item ``$r$'' is the integer dimension of the OpInf basis generated through POD. This is uniform across multiple subdomains.
    \item ``Data'' is the number of training steps taken, with a value of 1 only including the initial condition. There are a maximum of 101 steps possible across the time interval $[0,1]$, and any value of ``Data'' less than 101 results in prediction outside of training data. ``Data'' is uniform across all subdomains.
\end{enumerate}
Lastly, computational speed-up is measured as the wall clock time to run the online portion of the OpInf-Schwarz coupled ROM models relative to that of the fully FE coupled Schwarz method and the monolithic FE solution over the same domain configuration and on the same hardware. 

\subsection{Static Boundary Conditions}
\label{sec:results_simple}
For a first test, we consider the case where $u(\partial\Omega_T,t) = u(\partial\Omega_B,t) = 0$, $u(\partial\Omega_L,t) = 2$ and $u(\partial\Omega_R,t) = 5 \,\, \forall \, t \in (0,1]$. The discontinuity at some of the corners is resolved in favor of the top and bottom boundaries. We establish a baseline in Table~\ref{tab:simple_BCs_baseline} by coupling solely FE models in each of the 3 configurations of Figure~\ref{fig:schematics}. The wall time, in seconds, was averaged over 5 separate runs. Each subdomain has an overlap of 10 with its non-diagonal neighbors. The full monolithic simulation, averaged over 5 runs, takes 0.94 seconds to run the online portion.
\begin{table}[h]
    \begin{subtable}[h]{0.3\textwidth}
        \centering
        \begin{tabular}{c|c}
        $\mathbf{E}_{\ell^2}^{\text{avg}}$ &  $1.27 \times 10 ^{-14}$ \\ \hline Avg S.I. &  4.41 \\ \hline Time (s)  & 3.03
    \end{tabular}
          \caption{Vertical, Fig~\ref{schematic:vertical}}
     \label{tab:FE-FE_vertical_baseline}
    \end{subtable}
    \hfill
       \begin{subtable}[h]{0.3\textwidth}
        \centering
        \begin{tabular}{c|c}
        $\mathbf{E}_{\ell^2}^{\text{avg}}$ & $1.57\times 10 ^{-14}$ \\\hline Avg S.I. & 4.32 \\\hline Time (s) & 3.59 \\ 
         
    \end{tabular}
          \caption{Horizontal, Fig~\ref{schematic:horizontal}}
     \label{tab:FE-FE_hz_baseline}
    \end{subtable}
    \hfill
           \begin{subtable}[h]{0.3\textwidth}
        \centering
        \begin{tabular}{c|c}
        $\mathbf{E}_{\ell^2}^{\text{avg}}$ &  $1.94\times 10 ^{-14}$  \\ \hline Avg S.I. & 5.9 \\ \hline Time (s)  & 7.05
        
    \end{tabular}
          \caption{Four Squares, Fig~\ref{schematic:quadrants}}
     \label{tab:FE-FE_quads_baseline}
    \end{subtable}
     \caption{Baseline for all-FE coupled models, Overlap = 10}
     \label{tab:simple_BCs_baseline}
\end{table}

\begin{figure}[h]
     \centering
     \begin{subfigure}[h]{0.45\textwidth}
         \centering
         \includegraphics[width=\textwidth]{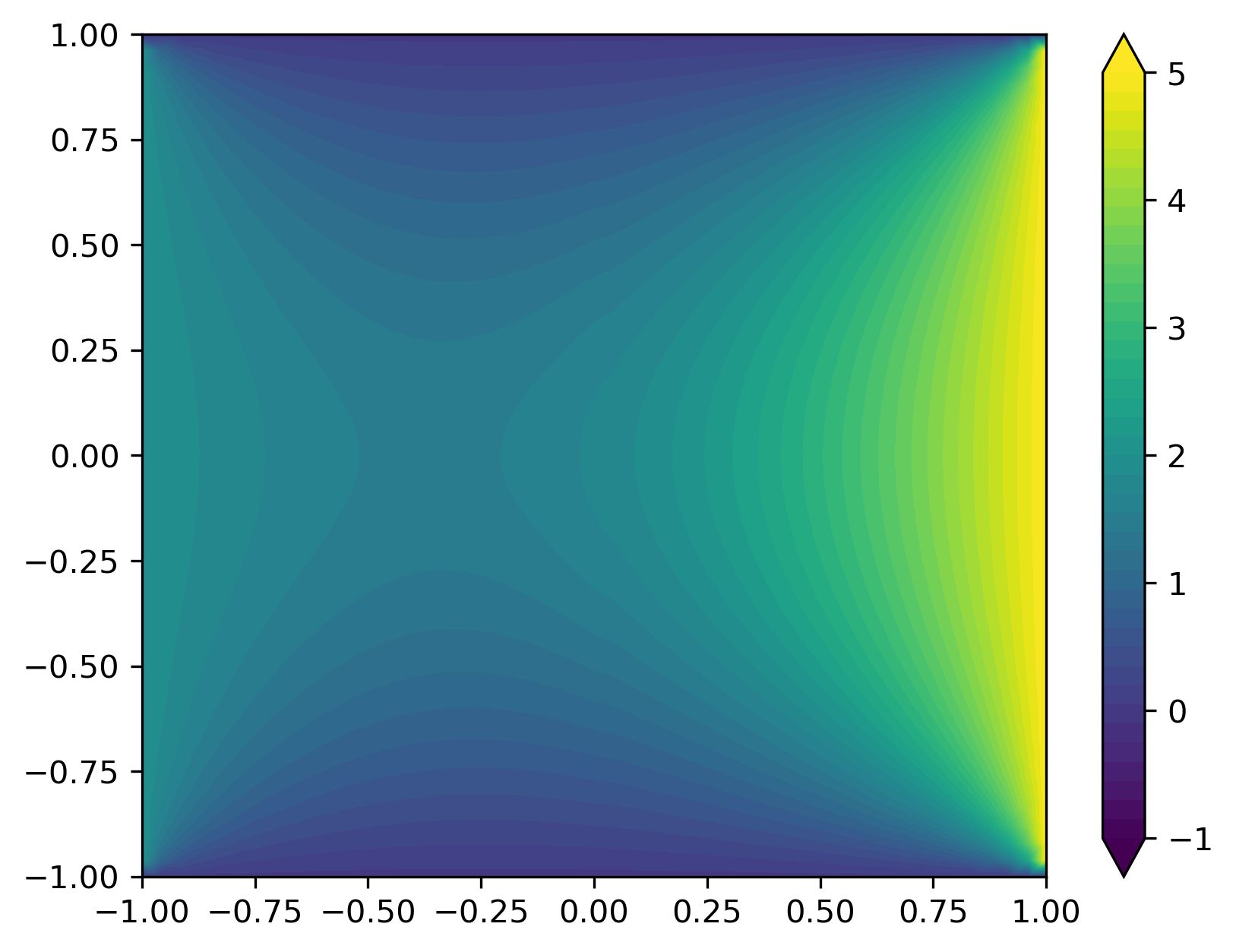}
         \caption{2 coupled OpInf models in vertical configuration, $t = 0$}
         \label{fig:simple_opinf}
     \end{subfigure}
     \hspace{5mm}
     \begin{subfigure}[h]{0.45\textwidth}
         \centering
         \includegraphics[width=\textwidth]{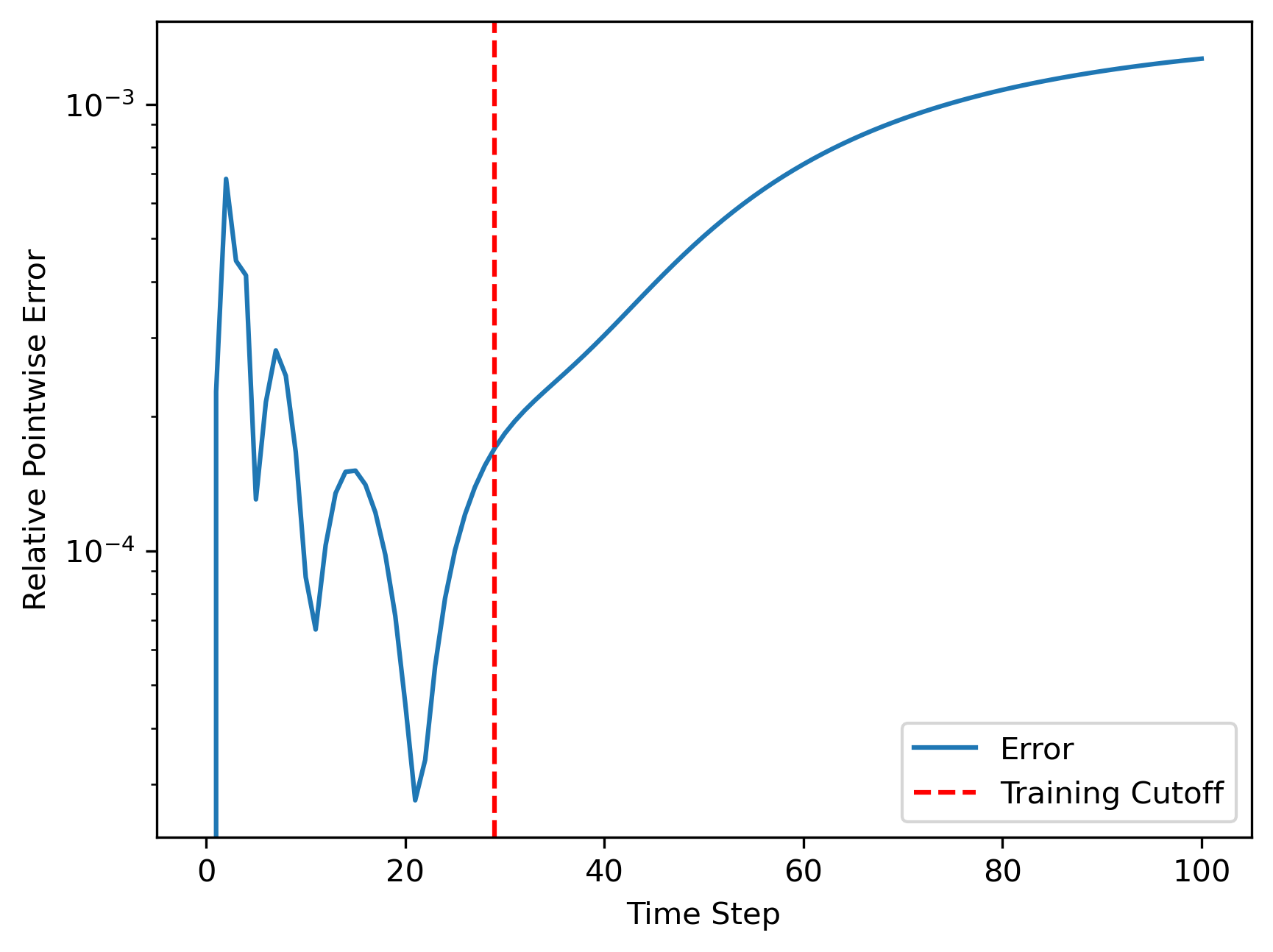}
         \caption{Error between OpInf-Schwarz and monolithic FE simulation over time.}
         \label{fig:simple_FOM_error}
     \end{subfigure}
\caption{OpInf-Schwarz vs. monolithic simulation. OpInf parameters: $r = 6$, Overlap = 10, Data = 30.}
\label{fig:simple_heat}
 \end{figure}

\paragraph{Vertical Configuration}
In Figure~\ref{fig:simple_heat} and Table~\ref{tab:OP-OP_vert}, we examine the behavior of the OpInf-Schwarz method with respect to $r$, the subdomain overlap, and data by coupling two OpInf models in the vertical orientation of Figure~\ref{schematic:vertical}.  Wall clock time results, in seconds, are averaged over 2 individual runs.

\begin{table}[h]
    \begin{subtable}[h]{0.45\textwidth}
        \centering
        \begin{tabular}{c | c| c|c}
Data & $\mathbf{E}_{\ell^2}^{\text{avg}}$ & Avg S.I. & Time \\ \hline
10 &$1.16\times 10^{-01}$ & 4.0 & 0.30 \\ \hline
20 &$5.41\times 10^{-03}$ & 4.57 & 0.30 \\ \hline
30 &$ 6.02\times 10^{-04}$ & 4.11 & 0.27\\ \hline
50 &$1.56\times 10^{-04}$ & 4.52 & 0.31 \\ \hline
80 & $9.16\times 10^{-05}$ & 4.39 & 0.28 \\ \hline
100 &  $8.62\times 10^{-05}$ & 4.49 & 0.31
       \end{tabular}
        \caption{Fixed parameters: Overlap = 10 a $r = 6$. }
        \label{tab:OP-OP_data_scaling}
     \end{subtable}
    \hspace{8mm}
    \begin{subtable}[h]{0.45\textwidth}
        \centering
        \begin{tabular}{c | c| c|c}
Overlap & $\mathbf{E}_{\ell^2}^{\text{avg}}$ & Avg S.I. & Time \\ \hline
1 &$8.68\times 10^{-03}$ & 6.4 & 0.35 \\ \hline
2 &$5.33\times 10^{-03}$ & 6.34 & 0.34 \\ \hline
5 &$ 2.36\times 10^{-03}$ & 5.7 & 0.31 \\ \hline
10 &$6.02\times 10^{-04}$ & 4.11 & 0.22 \\ \hline
20 & $5.96\times 10^{-04}$ & 4.7 & 0.25 \\ \hline
40 &  $8.18\times 10^{-04}$ & 4.26 & 0.24
       \end{tabular}
        \caption{Fixed parameters: Data = 30, $r = 6$. }
        \label{tab:OP-OP_overlap}
     \end{subtable}
  \\[\baselineskip]
    \centering
       \begin{subtable}[h]{0.6\textwidth}
        \centering
        \begin{tabular}{c | c| c|c|c}
$r$ & $\mathbf{E}_{\ell^2}^{\text{avg}}$ & Avg S.I. & $\mathbf{E}_{\text{proj}}^{\text{avg}}$ &Time \\ \hline
2 &$2.34\times 10^{-02}$ & 8.04 & $8.98\times 10^{-02}$ &  0.31 \\ \hline
4 &$4.24\times 10^{-03}$ & 5.0 & $1.02\times 10^{-02}$ & 0.25 \\ \hline
6 &$6.02\times 10^{-04}$ & 4.11 & $6.14\times 10^{-04}$ &0.19 \\ \hline
8 & $1.74\times 10^{-03}$ & 5.0 & $5.85\times 10^{-05}$ & 0.24 \\ \hline
10 & $1.85\times 10^{-03}$ & 5.0 & $6.02\times 10^{-06}$ & 0.23 \\ \hline
30 & $1.85\times 10^{-03}$ & 5.0 & $3.54\times 10^{-10}$ & 0.25
       \end{tabular}
       \caption{Fixed parameters: Data = 30, Overlap = 10.}
       \label{tab:OP-OP_r_scaling}
    \end{subtable}
    \caption{Comparison of OpInf-OpInf coupled models in vertical orientation (Fig~\ref{schematic:vertical}).}
     \label{tab:OP-OP_vert}
\end{table}


From a ROM perspective, Table~\ref{tab:OP-OP_r_scaling} has several interesting results which deserve explanation. First, $\mathbf{E}_{\ell^2}^{\text{avg}}$ only decreases up to $r = 6$, then increases before stalling out for higher values of $r$. This is despite the fact that the projection error of the snapshots onto the basis drops consistently with the basis size. A basis size of  $r = 3$ and $r = 4$ is necessary to capture 99.9\% of the total energy of the snapshots $\Omega_1$ and $\Omega_2$ respectively. 

In the row $r = 6$, $\mathbf{E}_{\ell^2}^{\text{avg}}$ and $\mathbf{E}_{\text{proj}}^{\text{avg}}$ have the same order of accuracy, while $\mathbf{E}_{\text{proj}}^{\text{avg}}$ is lower for all further values of $r$. This suggests that the accuracy of the OpInf-Schwarz method is limited by the basis quality up until $r = 6$, and afterwards is most limited by error in the inferred operators. This could be due to several reasons, such as data quality, our regularization strategy affecting accuracy, or perhaps the Schwarz coupling method interacting with the ROM in an unexpected way. 

Another item of note in Table~\ref{tab:OP-OP_r_scaling} is that the average run time decreases by about a factor of three as we increase $r$ from 2 to 6, even though we are increasing the size of the linear systems we are solving. This is due to the coupling of poor quality solutions when $r=2$, which requires additional Schwarz iterations to converge as shown by the mean Schwarz iterations.  A similar phenomenon is observed in~\cite{Barnett:2022Schwarz} when coupling subdomain-local intrusive ROMs via the Schwarz alternating method.


The remaining subtables of Table~\ref{tab:OP-OP_vert} reveal that OpInf-Schwarz performs as anticipated. In Table~\ref{tab:OP-OP_overlap}, increasing the size of the overlap decreases the number of Schwarz iterations, and in Table~\ref{tab:OP-OP_data_scaling}, increasing the training data decreases the approximation error. Since we are using the same data for both the POD basis and the OpInf training, it is fair to ask whether increasing the data allowance also allows POD to extract useful data for larger values of $r$---that is, further decreases to $\mathbf{E}_{\text{proj}}^{\text{avg}}$ . In our testing, increasing training data had a minor effect on $\mathbf{E}_{\text{proj}}^{\text{avg}}$ compared to increasing $r$, though this would likely become more relevant in more challenging problems. The minimal error we were able to recover was $\mathbf{E}_{\ell^2}^{\text{avg}} = 2.5 \times 10^{-05}$ with the choices of $r = 10$, Overlap = 10, and 100 steps of training data.  

While the times recorded are certainly not deterministic, a comparison between the FE-FE coupling of Table~\ref{tab:FE-FE_vertical_baseline} and the minimum and maximum values of Table~\ref{tab:OP-OP_vert} suggest a speed-up of 8.6 to 15.2 times for fully OpInf-OpInf coupled models as compared to a FOM-FOM coupling on the same spatio-temporal grid, and 2.7 to 4.9 times speed-up for fully OpInf-OpInf coupled models as compared to the monolithic FE simulation.  Further speed-ups may be possible by introducing the ``additive'' Schwarz method~\cite{wentland2024Schwarz}, which computes subdomain solutions in parallel.   


\paragraph{Comparison With Other Configurations}
We now investigate whether interactions between the subdomain geometry and boundary conditions have an effect on the performance of OpInf-Schwarz. Because the non-zero boundary conditions in this simple case are on the left and right, a vertical overlap splitting the middle has a less complicated gradient than a horizontal overlap. In the scenario of $r = 6$, Overlap = 10, and 30 steps of training data, Table~\ref{tab:simple_vertical_square} shows that the error and Schwarz iterations are mildly worse for those configurations that impose a more complicated Schwarz boundary. 

\begin{table}[h]
    \centering
    \begin{tabular}{c|c|c|c}
        Model Configuration & $\mathbf{E}_{\ell^2}^{\text{avg}}$ & Avg S.I. & Time (s)  \\ \hline
        Vertical & 6.02 $\times 10^{-04}$  & 4.11 & 0.26 \\ \hline
        Horizontal & 8.91 $\times 10^{-04}$  & 4.17 & 0.26 \\ \hline
        Four Squares &7.10 $\times 10^{-04}$  & 5.81 & 0.61
    \end{tabular}
    \caption{Vertical, horizontal, and square configurations for simple boundaries}
    \label{tab:simple_vertical_square}
\end{table}

All of the results in the simple case of Section~\ref{sec:results_simple} are in line with prior expectations for Schwarz performance in a ROM setting. 

\subsection{Time-Varying Boundary Conditions}
\label{sec:results_time_var}
We now consider a more challenging scenario for the 2D heat equation. Instead of solely constant boundary conditions, define 
\begin{equation}
\label{eq: bc_time_var}
    q(t, \mu) = 1 + 0.5 \sin(2 \pi\mu t),
\end{equation}
to be a time varying boundary condition parameterized by the frequency $\mu$. We consider the monolithic heat equation of $\eqref{eq:heat_pde}$ on $\Omega$ specified by $u(\partial\Omega_L,t) = q(t,2)$, $u(\partial\Omega_T,t) = q(t,4)$, $u(\partial\Omega_B,t) = 5$, $u(\partial\Omega_R,t) = 1 \,\, \forall t \in (0,1]$ and $u(\Omega,0) = 0$.

The offset frequencies of the boundary conditions on the left and top boundaries create dynamics which evolve through a significant portion of the time domain, especially in the top left quadrant. These types of problems can be challenging for ROMs to resolve when given insufficient data, though they can usually be fixed by increasing the training data. This is one of the essential choices in ROMs -- increasing the training range can often greatly improve accuracy, but data generation in realistic problems can be enormously expensive, especially if it must be generated from a monolithic simulation. Domain decomposition frameworks provide options beyond increasing the training range, and we now consider a comparison between purely OpInf-Schwarz with a FE-OpInf-Schwarz coupled problem. 

\begin{figure}[h]
     \centering
     \begin{subfigure}[h]{0.45\textwidth}
         \centering
         \includegraphics[width=\textwidth]{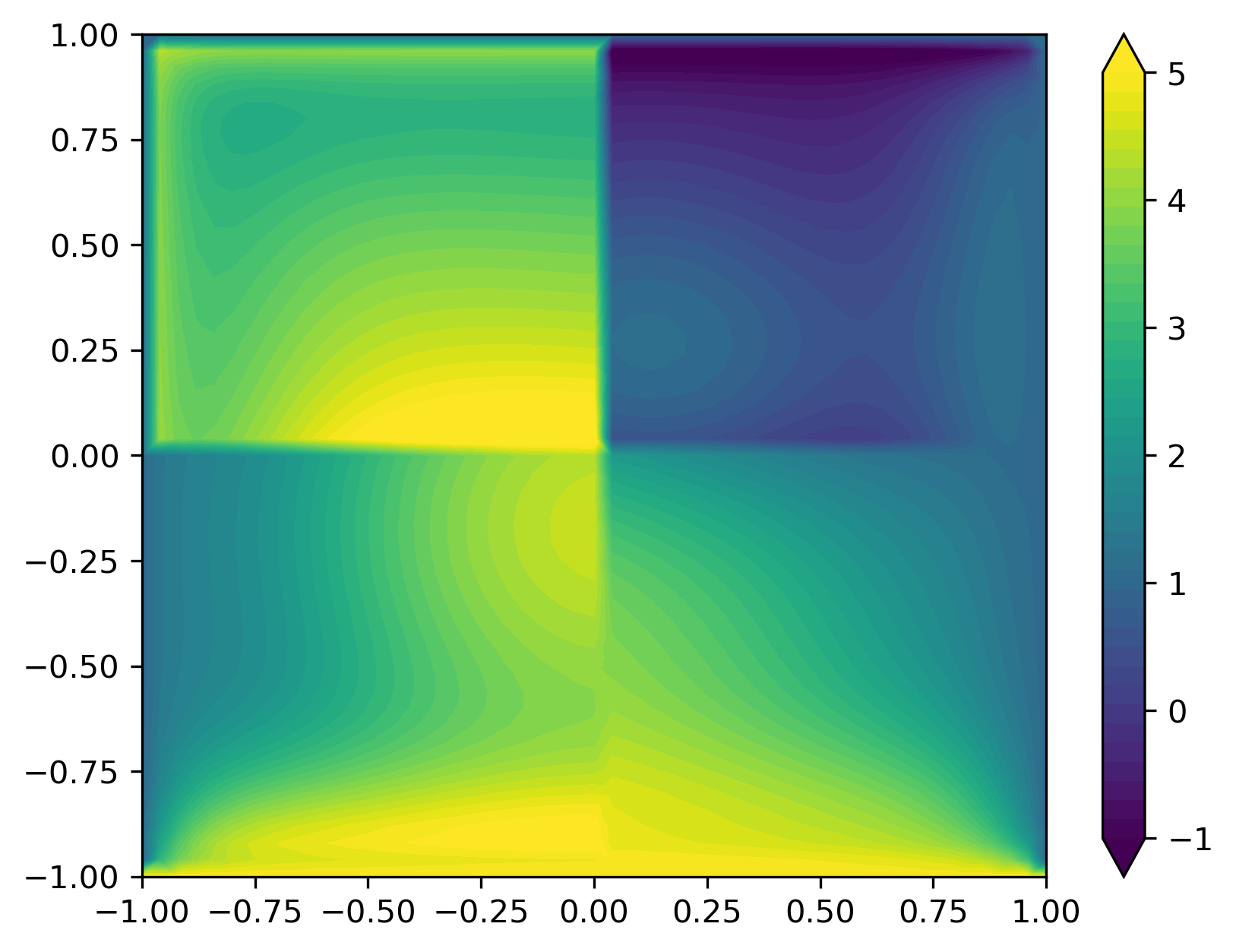}
         \caption{OpInf solution in all 4 subdomains, 40 steps of training data in all subdomains.}
         \label{fig:time_varying_opinf_all}
     \end{subfigure}
     \hspace{5mm}
     \begin{subfigure}[h]{0.45\textwidth}
         \centering
         \includegraphics[width=\textwidth]{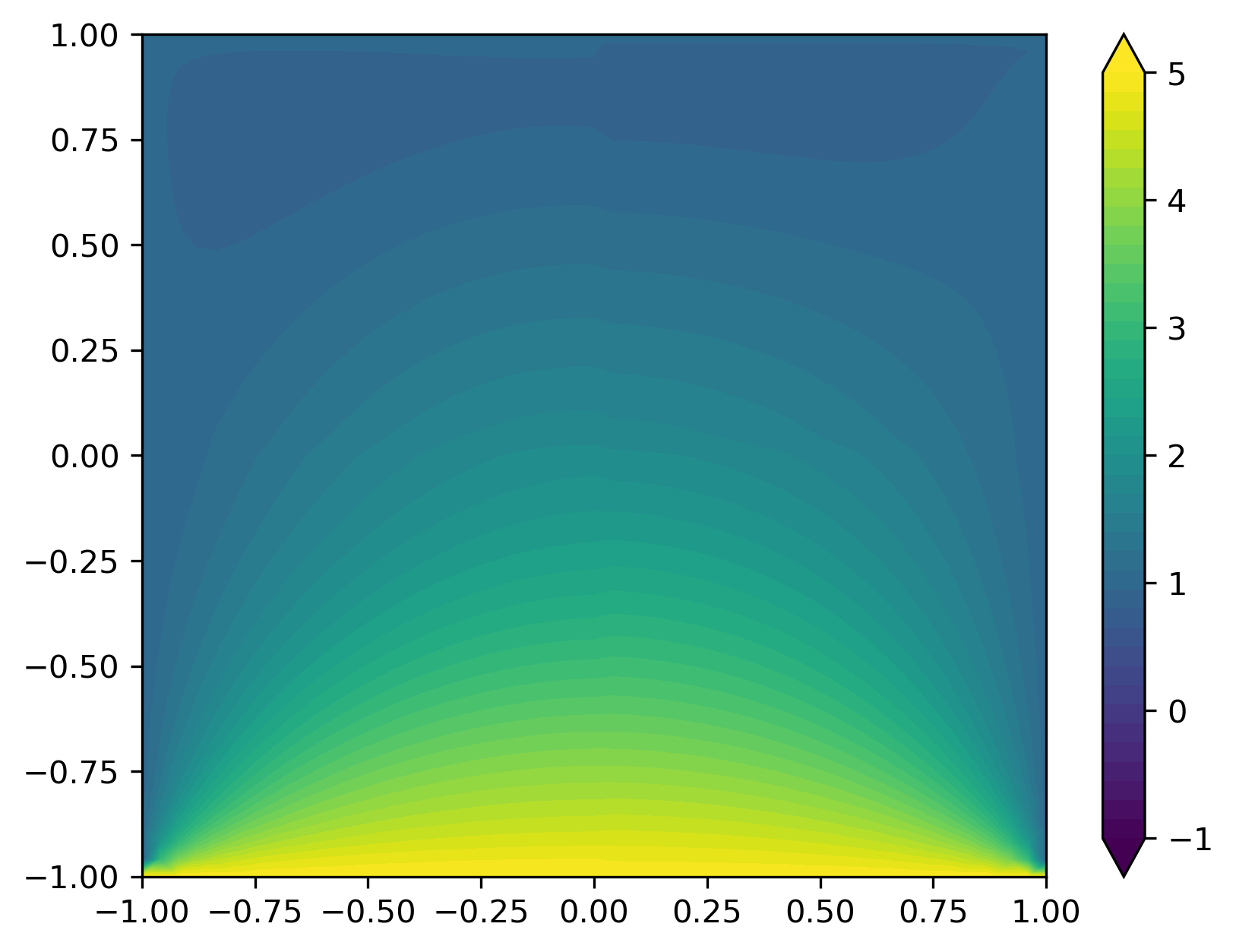}
         \caption{FE solution in top left $\Omega_3$ only, 40 steps of training data for OpInf portions.}
         \label{fig:time_varying_fom_top_left}
     \end{subfigure}
\\
     \begin{subfigure}[h]{0.45\textwidth}
         \centering
         \includegraphics[width=\textwidth]{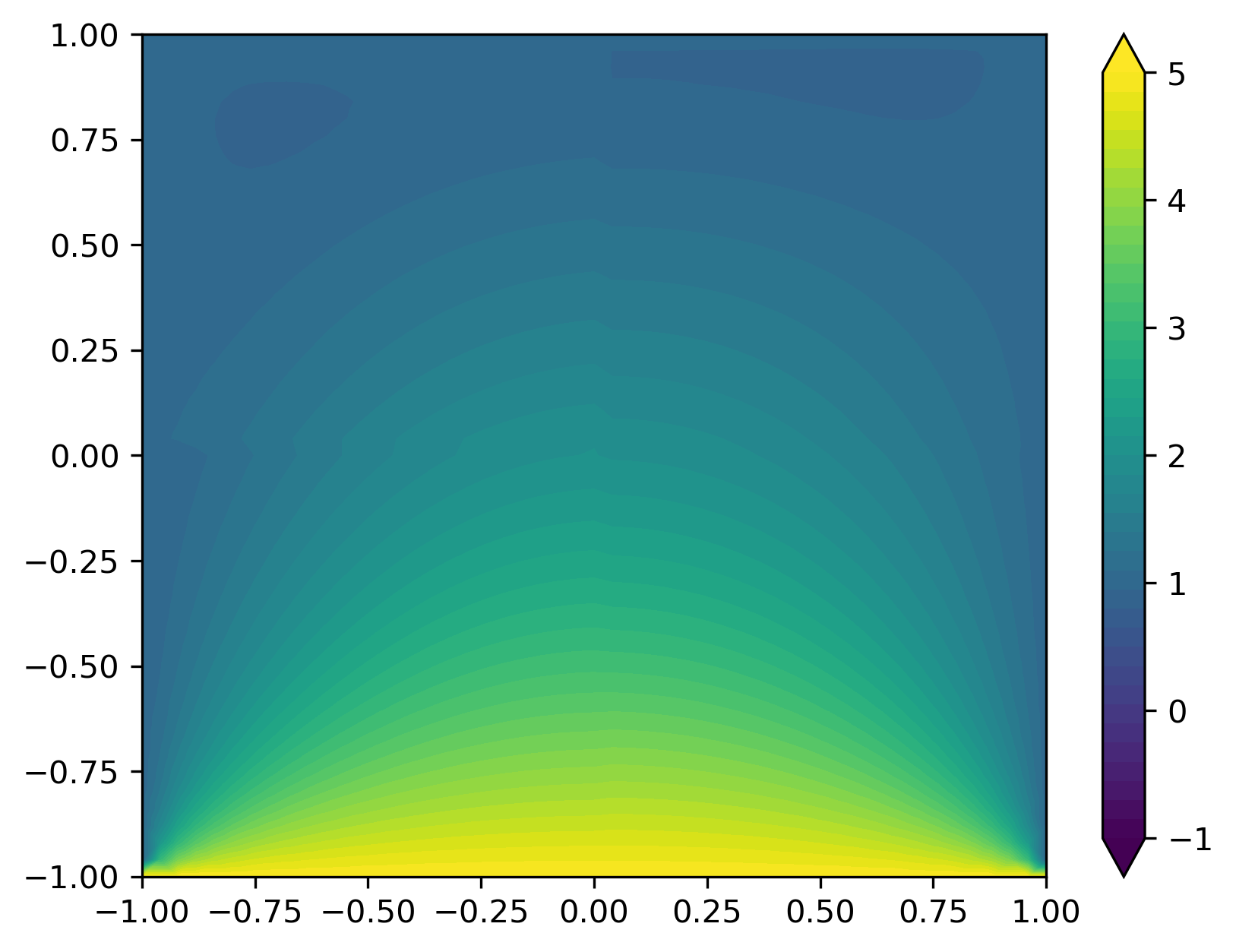}
         \caption{OpInf solution in all 4 subdomains, 65 steps of training data in all subdomains.}
         \label{fig:time_varying_opinf_all_65_data}
     \end{subfigure}
     \hspace{5mm}
     \begin{subfigure}[h]{0.45\textwidth}
         \centering
         \includegraphics[width=\textwidth]{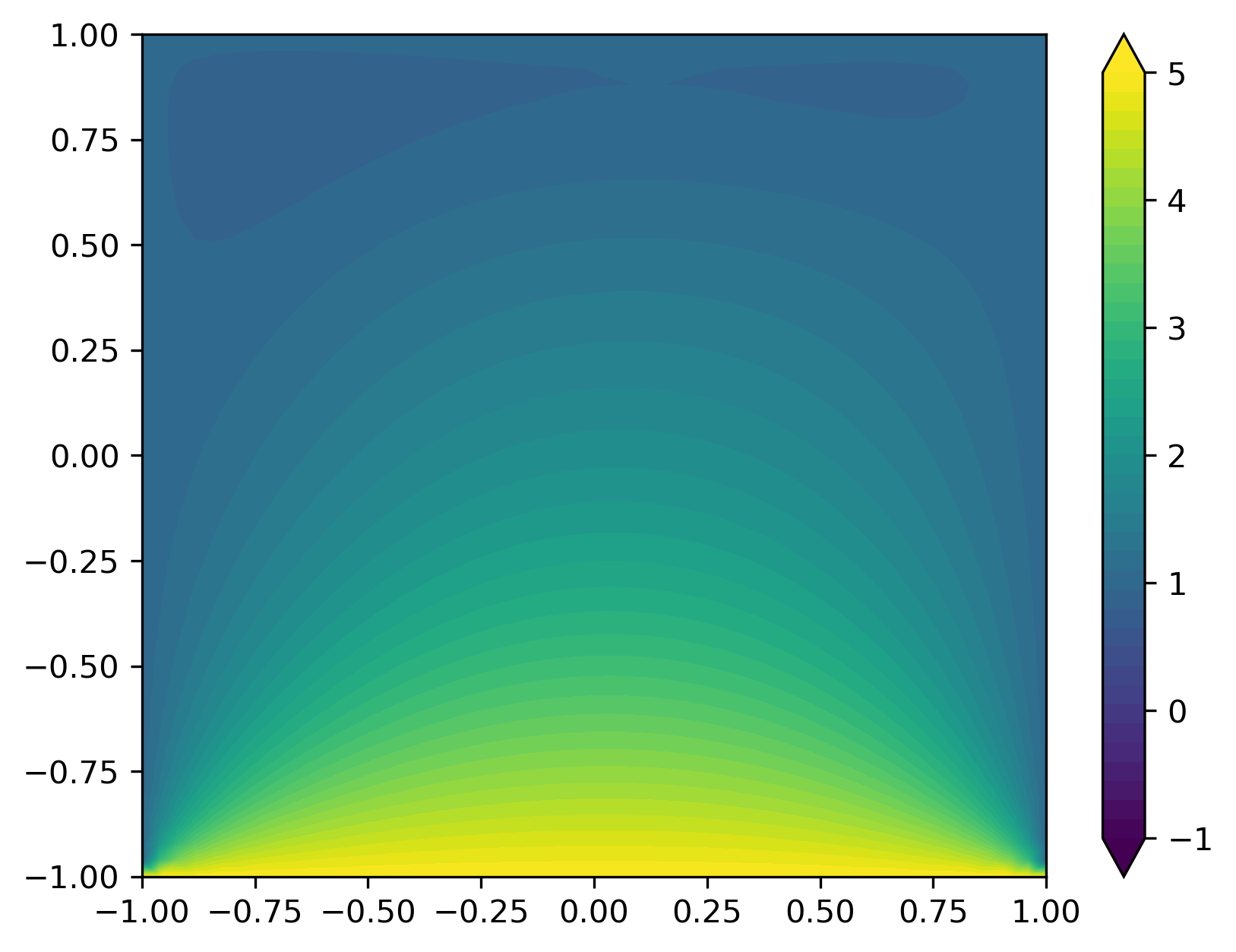}
         \caption{Monolithic FOM (Finite Element Simulation).}
         \label{fig:time_varying_monolithic}
     \end{subfigure}
     
        \caption{OpInf-Schwarz solutions and Monolithic FOM at $t = 1$.}
        \label{fig:time_varying_schwarz}
\end{figure}

\paragraph{Four Squares Configuration} Figure~\ref{fig:time_varying_schwarz} shows two OpInf-Schwarz solutions at the final time step of $t = 1$ in the four square configuration of Figure~\ref{schematic:quadrants}. For the OpInf portions, each solution has the same ROM basis dimension $r = 6$, 40 steps of training data (out of a possible 101), and a mean centering strategy performed on the snapshots. They additionally share the same boundary and initial conditions as described above. The only difference in their setup is that Figure~\ref{fig:time_varying_opinf_all}
 couples four OpInf models, while Figure~\ref{fig:time_varying_fom_top_left} has replaced one of the OpInf models with a FE model in $\Omega_3$, the top-left subdomain. 

 
 The FE assimilation in Figure~\ref{fig:time_varying_fom_top_left} is not perfect if one looks closely along $x = 0$ or $y = 0$, but it shows a cohesive solution very similar to the monolithic solution in Figure~\ref{fig:time_varying_monolithic}, while the solution in Figure~\ref{fig:time_varying_opinf_all} is clearly failing and is in fact about to diverge entirely. The All OpInf model with a higher data allowance in Figure~\ref{fig:time_varying_opinf_all_65_data} is comparable to Figure~\ref{fig:time_varying_fom_top_left}. 

 One interesting feature of Figure~\ref{fig:time_varying_opinf_all} is that the Schwarz communication pattern is clearly evident.  The solution is completely incorrect in the top left subdomain, $\Omega_3$.  The error in $\Omega_3$ subsequently propagates to neighboring $\Omega_1$ and $\Omega_4$, while the bottom right subdomain $\Omega_2$ remains closest to an accurate solution. 

\begin{table}[h]
    \begin{subtable}[h]{0.4\textwidth}
        \centering
        \begin{tabular}{c | c| c}
Data& $\mathbf{E}_{\ell^2}^{\text{avg}}$ &  $\mathbf{E}_{\ell^2}^{\text{max}}$   \\ \hline
30 &$1.97\times 10^{-01}$ & $1.22\times 10^{00}$  \\ \hline
35 &$5.41\times 10^{-01}$ & 4.07$\times 10^{-01}$  \\ \hline
40 &$1.07\times 10^{-01}$ & $6.40\times 10^{-01}$  \\ \hline
50 & $3.62\times 10^{-02}$ &$2.12\times 10^{-01}$  \\ \hline
60 & $9.15\times 10^{-03}$ &  $3.13\times 10^{-02}$  \\ \hline
65 & $6.48\times 10^{-03}$ & $1.91\times 10^{-02}$ 
       \end{tabular}
       \caption{All OpInf, 4 Coupled OpInf Models}
       \label{tab:squares_OP_all}
    \end{subtable}
    \hfill
    \begin{subtable}[h]{0.45\textwidth}
        \centering
        \begin{tabular}{c | c| c}
Data & $\mathbf{E}_{\ell^2}^{\text{avg}}$ & $\mathbf{E}_{\ell^2}^{\text{max}}$  \\ \hline
30 &$3.33\times 10^{-02}$ & $7.06\times 10^{-02}$ \\ \hline
35 &$9.47\times 10^{-03}$ & $2.75\times 10^{-02}$  \\ \hline
40 &$ 7.14\times 10^{-03}$ & $1.89\times 10^{-02}$  \\ \hline
50 &$4.76\times 10^{-03}$ & $1.11\times 10^{-02}$ \\ \hline
60 & $1.86\times 10^{-03}$ & $1.4\times 10^{-02}$ \\ \hline
65 &  $1.77\times 10^{-03}$ & $1.40\times 10^{-02}$ 
       \end{tabular}
        \caption{3 OpInf 1 FE, FE model in top left subdomain $\Omega_3$ only.}
        \label{tab:squares_OP_FOM}
     \end{subtable}
     \hfill
    
     \caption{OpInf-Schwarz errors for various training data ranges with $r = 6$ and Overlap = 10.}
     \label{tab:Time_Varying_Squares}
\end{table}

From Table~\ref{tab:Time_Varying_Squares}, it can be inferred that a purely OpInf coupled Schwarz model would require 65 total timesteps (an additional 25 steps of training data) to exceed the average relative accuracy of the OpInf-Schwarz with one FE subdomain and 40 steps of training data (as is displayed in Figure~\ref{fig:time_varying_fom_top_left}).
However, over two separate runs, the average time taken for the 4 coupled OpInf models in Table~\ref{tab:squares_OP_all} was 0.73 seconds, while the average for the substitution of one FE model in $\Omega_3$ in Table~\ref{tab:squares_OP_FOM} is 3.15 seconds, a 4.3 times slowdown over a purely OpInf coupled model. The relative ease of model switching in the Schwarz framework lets one choose between data generation costs and monolithic FOM costs.

In comparison, a completely FE coupled model with the same settings has $\mathbf{E}_{\ell^2}^{\text{avg}} = \mathcal{O}(10^{-14})$, but, over three runs, takes an average of 7.9 seconds. A monolithic FE sim takes approximately 0.8 seconds. The ``All OpInf" model is therefore about 9 times faster, and the 3 OpInf 1 FE model about 2.4 times faster than the fully FE coupled model. 
Comparing against the monolithic FOM, the ``All OpInf" model is 1.1 times faster and the 3 OpInf 1 FE is 3.9 times slower. 

\paragraph{Efficiency Comparison With Other Configurations} 
\begin{figure}[h]
     \centering
 \includegraphics[width=\textwidth]{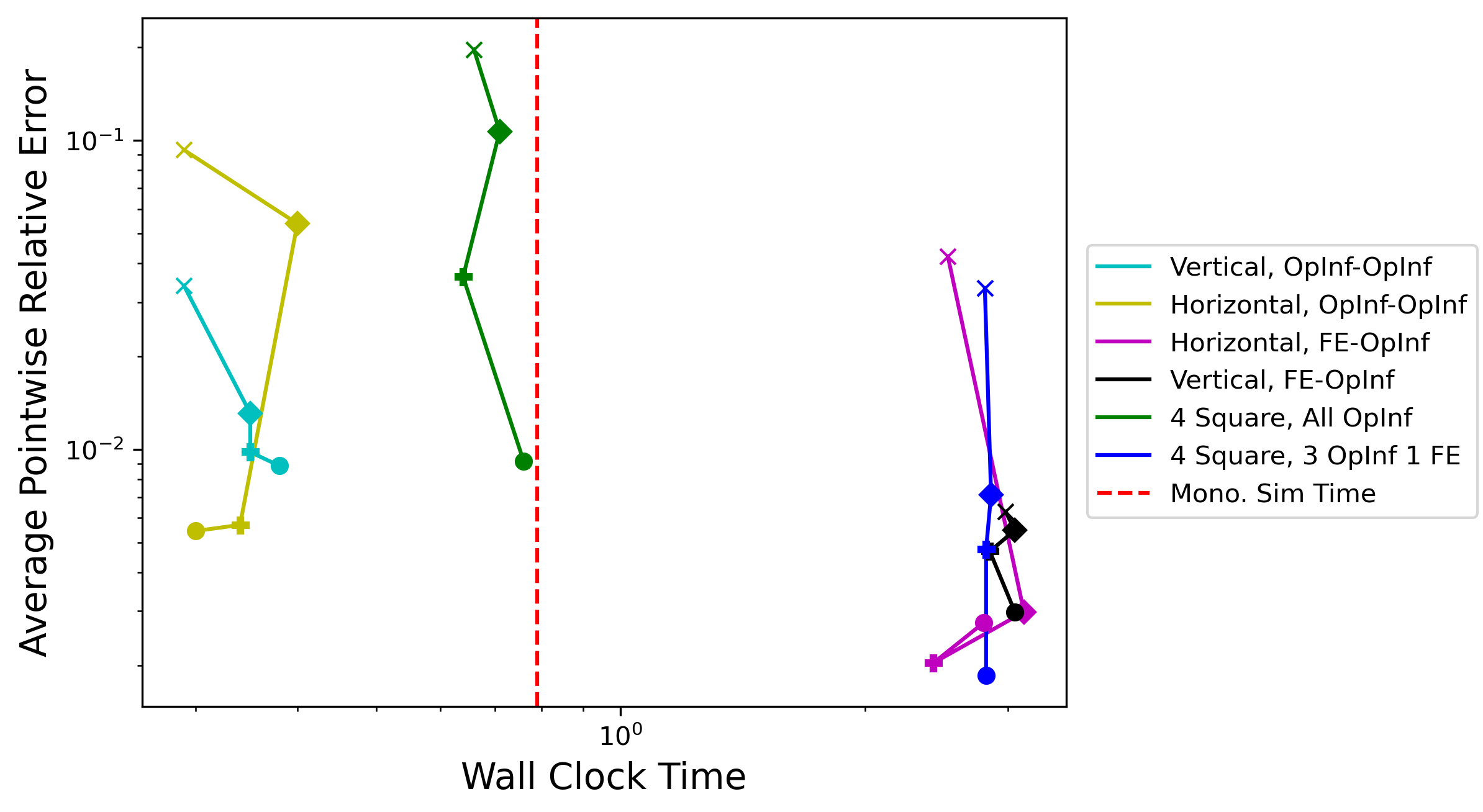}
\caption{Pareto plot for time-varying BCs, $\mathbf{E}_{\ell^2}^{\text{avg}}$ vs Time (s). Fixed parameters: $r = 6$, Overlap $=10$. Plotted points are ``Data'' as it ranges over 30, 40, 50, and 60 marked by $\times$, $\blacklozenge$, $\boldsymbol{+}$, and $\bullet$ respectively as in Table~\ref{tab:Time_Varying_Squares}. Values in bottom left are preferred.}
    \label{fig:pareto}
 \end{figure}

We summarize our timing results for the time-varying boundary problem in the Pareto plot of Figure~\ref{fig:pareto} using data from Table~\ref{tab:Time_Varying_Squares} for the ``Four Squares'' labelled results, and additional results for the ``Vertical'' labelled results produced in the vertical configuration from Figure~\ref{schematic:vertical}. In particular, ``Vertical FE-OpInf'' was produced with the FE model in the left subdomain, and ``Horizontal FE-OpInf'' was produced with the FE subdomain as the top subdomain, both advantageous DDs for resolving the more complicated dynamics in the top left portion of $\Omega$. Fixed parameters are $r = 6$ and Overlap = 10. This chart demonstrates the significant impact that subdomain geometry can have when interacting with more complicated boundary issues or solution features. The vertical OpInf-OpInf is more accurate for smaller data allowances, but it is overtaken in accuracy by the 3 OpInf 1 FE model for higher data allowances. All of the models with a FE component have similar time costs, but the horizontal and four square configurations make better use of additional data than the vertical configuration does. 

While none of the models incorporating a FE subdomain exhibit lower runtimes than that of the monolithic FOM, we reiterate that it is possible to improve the performance of the Schwarz iteration process via the additive Schwarz formulation.

\section{Conclusions and Future Work} \label{sec:conclusions}

The major impetus driving the development of the OpInf-Schwarz method is to achieve a minimally intrusive methodology for constructing ROMs that can interface modularly with existing high performance codes. A truly data-driven and modular domain decomposition-based ROM has the potential to greatly reduce computational expenses incurred by meshing and re-meshing complicated 3D objects, and therefore in performing long run-time and/or multi-query simulations. 

The preliminary investigation presented herein has assessed various characteristics of the OpInf-Schwarz framework in the context of a 2D heat equation. Our results have shown that the OpInf-Schwarz method is capable of coupling  OpInf and FOM subdomains to recover accurate solutions, and that a purely OpInf formulation can run faster than a monolithic simulation.  We have also demonstrated that OpInf-Schwarz is very flexible, and can be tuned in several ways to the problem at hand (e.g., by choosing a different DD or ROM/FOM assignment to the subdomains).

At the same time, this work has uncovered several challenges that are worth investigating in the future. First, the calculated average relative pointwise error presented in Section~\ref{sec:results_simple} stagnates at around $\mathcal{O}(10^{-5})$ and does not decrease further despite a increase to the POD basis size and training data. While this level of error is acceptable for many applications, it would be ideal for the method to be truly convergent to the data source. While it seems likely that the issue lies with the error of the inferred operators from the  regression problem~\eqref{eq:opinf-schwarz_regression_problem} given that the projection errors reported in Table~\ref{tab:OP-OP_r_scaling} continue to decrease, a comparison with a standard intrusive projection ROM would confirm whether the cause is error in the inferred operators or errors induced by the coupling process. It may be possible to reduce optimization errors by employing 
a technique called re-projection~\cite{Peherstorfer:2020}.

There are also potential improvements to the implementation of OpInf presented in this paper, as the stability of the learned operators is not enforced at the inference step. At present, this is resolved with a simple Tikhonov-based regularization strategy, but a more rigorous approach may be needed in more complicated problems, e.g., by using a more sophisticated regularization strategy such as those presented in~\cite{McQuarrie2021combustion, Sawant:2023}. There are additional opportunities to improve the implementation of the Schwarz boundary conditions within our OpInf-Schwarz formulation in a way that both reduces computational cost and improves accuracy.  

Currently, the size of the learned boundary operator  $\hat{\mathbf{B}}$ scales with the size of the boundary in our formulation, but considering that even in this simple 2D case the number of boundary nodes is already 20 or more times larger than the ROM dimension, it is worth questioning whether all of this boundary information is being used in a meaningful way. 
Reducing the computational complexity of our method requires reducing the number of columns of $\hat{\mathbf{B}}$, which can be large for multi-dimensional problems with many boundary nodes.  To mitigate the cost associated with evaluating $\hat{\mathbf{B}} \mathbf{g}$, it may be possible to perform a separate dimension reduction of this term.  Finally, as noted earlier, the overall computational complexity of OpInf-Schwarz can be improved by introducing ``additive'' Schwarz, characterized by parallel subdomain solves with asynchronous boundary condition communication~\cite{wentland2024Schwarz}. 


The last major point is on extending the Opinf-Schwarz method to more challenging and realistic problems, as we expect the greatest utility and time savings of this model to be shown in more complicated situations. The heat equation is a standard prototype for model reduction methods, but extending the model to 3D, non-linear, convection dominated, and parametric problems~\cite{McQuarrie2023parametric,Vijaywargiya:2024,yildiz2021parametric} is a priority, as are implementing strategies to choose optimal domain decomposition and model assignment to enable on-the-fly ROM-FOM switching for maximal performance.

\section*{Acknowledgements}

This material is based upon work supported by the U.S. Department of Energy, Office of Science, Office of Advanced Scientific Computing Research, Mathematical Multifaceted Integrated Capability Centers (MMICCs) program, under Field Work Proposal 22-025291 (Multifaceted Mathematics for Predictive Digital Twins (M2dt)), Field Work Proposal 20-020467, and the Laboratory Directed Research and Development program at Sandia National Laboratories. The writing of this manuscript was funded in part by the fourth author’s (Irina Tezaur’s) Presidential Early Career Award for Scientists and Engineers (PECASE). This article has been authored by an employee of National Technology \& Engineering Solutions of Sandia, LLC under Contract No. DE-NA0003525 with the U.S. Department of Energy (DOE). The employee owns all right, title and interest in and to the article and is solely responsible for its contents. The United States Government retains and the publisher, by accepting the article for publication, acknowledges that the United States Government retains a non-exclusive, paid-up, irrevocable, world-wide license to publish or reproduce the published form of this article or allow others to do so, for United States Government purposes. The DOE will provide public access to these results of federally sponsored research in accordance with the DOE Public Access Plan \url{https://www.energy.gov/downloads/doe-public-access-plan}.

The authors would like to thank Dr.~Shane McQuarrie for numerous helpful technical discussions in the area of operator inference, and for providing support on the {\tt OpInf} library.

\bibliographystyle{siam}
\bibliography{IanMoore}

\end{document}